\documentclass[10pt]{article}
\usepackage{amsthm}
\usepackage{fullpage}
\usepackage{amssymb, amsmath}
\usepackage{hyperref}
\usepackage{graphicx}
\usepackage{titlesec}
\usepackage{enumitem}
\usepackage{comment}
\usepackage{relsize}
\usepackage[font=small]{caption}
\newtheorem{thm}{Theorem}
\newtheorem{prop}[thm]{Proposition}
\newtheorem{lem}[thm]{Lemma}

\newtheorem{conj}[thm]{Conjecture}
\newtheorem{ques}{Question}

\theoremstyle{definition}
\newtheorem{example}[thm]{Example}

\title{Bounded degree graphs and hypergraphs with no full rainbow matchings}
\author{Ronen Wdowinski\thanks{Department of Combinatorics and Optimization, University of Waterloo, Waterloo ON, Canada. Email: ronen.wdowinski@uwaterloo.ca} \thanks{Institute of Discrete Mathematics, Graz University of Technology, Graz, Austria. Email: wdowinski@math.tugraz.at}}  
\date{\today}

\begin{document}

\maketitle

\begin{abstract}
Given a multi-hypergraph $G$ that is edge-colored into color classes $E_1, \ldots, E_n$, a full rainbow matching is a matching of $G$ that contains exactly one edge from each color class $E_i$. One way to guarantee the existence of a full rainbow matching is to have the size of each color class $E_i$ be sufficiently large compared to the maximum degree of $G$. In this paper, we apply an iterative method to construct edge-colored multi-hypergraphs with a given maximum degree, large color classes, and no full rainbow matchings. First, for every $r \ge 1$ and $\Delta \ge 2$, we construct edge-colored $r$-uniform multi-hypergraphs with maximum degree $\Delta$ such that each color class has size $|E_i| \ge r\Delta - 1$ and there is no full rainbow matching, which demonstrates that a theorem of Aharoni, Berger, and Meshulam (2005) is best possible. Second, we construct properly edge-colored multigraphs with no full rainbow matchings which disprove conjectures of Delcourt and Postle (2022). Finally, we apply results on full rainbow matchings to list edge-colorings and prove that a color degree generalization of Galvin's theorem (1995) does not hold.
\end{abstract}

\section{Introduction}
A multi-hypergraph $G$, which has vertex set $V(G)$ and edge multiset $E(G)$, is \textit{$r$-uniform} if every edge of $G$ has size $r$, and in this case we will refer to $G$ as an \textit{$r$-graph}. A multi-hypergraph $G$ is \textit{$k$-partite} if there exists a partition of $V(G)$ into $k$ parts $V_1, \ldots, V_k$, such that each edge intersects each part $V_i$ in at most one vertex. The \textit{degree} of a vertex $v$ is the number of edges of $G$ containing $v$. A \textit{matching} of $G$ is a set of pairwise disjoint edges. Suppose that $E_1, \ldots, E_n$ is the set of color classes of a (not necessarily proper) edge-coloring of $G$, meaning that $E_1, \ldots, E_n$ is any partition of $E(G)$. A \textit{rainbow matching} with respect to $E_1, \ldots, E_n$ is a matching of $G$ that contains at most one edge from each color class $E_i$. A rainbow matching is \textit{full} if it contains exactly one edge from each color class $E_i$. A full rainbow matching can also be viewed as a system of pairwise-disjoint representative edges from a family of hypergraphs $\mathcal{E}$ on the same vertex set.

Research on large or full rainbow matchings in edge-colored multi-hypergraphs has in large part been motivated by problems about transversals in Latin squares (see \cite{Wa} for a more in-depth overview). The famous problem of finding an orthogonal mate to a Latin square is equivalent to finding a decomposition of that Latin square into transversals. However, the question of whether there exists even a single large partial transversal in a given Latin square is already interesting. The celebrated conjecture of Ryser, Brualdi, and Stein \cite{BrRy, Ry, St} states that in fact every Latin square has a partial transversal that misses at most one symbol (a proof for large $n$ was recently announced by Montgomery \cite{Mo}). An $n \times n$ Latin square is well-known to be equivalent to a proper edge-coloring of $K_{n,n}$ with $n$ colors, and then a partial transversal of size $m$ in the Latin square corresponds to a rainbow matching of size $m$ in the edge-colored $K_{n,n}$. Hence, many related conjectures and theorems have sprung up about full rainbow matchings in general edge-colored multigraphs, as well as in multi-hypergraphs. Some well-known conjectures are those of Aharoni and Berger \cite{AhBe} and of Gao, Ramadurai, Wanless, and Wormald \cite{GaRaWaWo}, who assert the existence of a full rainbow matching in a properly edge-colored multigraph $G$ if there are $n$ color classes each of size at least $n+1$ (if $G$ is bipartite) or $n+2$ (if $G$ is not bipartite), respectively. These conjectures have received a lot of attention, with many partial and asymptotic results \cite{BaGySa, GaRaWaWo, HaSh, KePoSuYe, KeYe, Po} (see \cite{MuPoSu} for short proofs).

In this paper, we address a different but related problem: the existence of full rainbow matchings in edge-colored multi-hypergraphs when we only assume that the color class sizes are sufficiently large compared to the maximum degree of the multi-hypergraph, irrespective of the number of color classes. This angle on rainbow matchings has received some attention in recent years \cite{AhBe, AhBeKoZi, ChLo, DePo, GaRaWaWo}. Our main contributions in this area are various, often optimal, constructions of edge-colored multi-hypergraphs with a given maximum degree, large color class sizes, and no full rainbow matchings. These constructions come from a simple iterative method (based on Lemma \ref{join-lemma} below) that was developed and applied extensively by Haxell and this paper's author in \cite{HaWd1} (as well as in \cite{CaHaKaWd, HaWd2}), for the related topic of \textit{independent transversals} in graphs. Independent transversals can be viewed as ``full rainbow independent sets" in vertex-colored graphs, and a full rainbow matching in a multi-hypergraph is the same as an independent transversal in the line graph of that multi-hypergraph. Maximum degree conditions for the existence of independent transversals in graphs are prevalent in the literature (e.g. \cite{AhAlBe, AhHoHoSp, Ha2, HaSz1, LoSu, SzTa}), and the method developed in \cite{HaWd1} has yielded optimal constructions of vertex-partitioned graphs with no independent transversals. The method generalizes to all simplicial complexes, and in this paper we apply it to matching complexes of multi-hypergraphs.

First, we focus on full rainbow matchings in general edge-colored $r$-graphs with a given maximum degree $\Delta$. One sufficient condition for the existence of full rainbow matchings is the following theorem, which follows directly from a Hall-type theorem of Aharoni, Berger, and Meshulam \cite{AhBeMe} (see Section \ref{proof-of-thm-1} for details).

\begin{thm}[Aharoni, Berger, Meshulam] \label{max-degree-thm}
Let $G$ be an $r$-graph with maximum degree $\Delta$, and let $E_1, \ldots, E_n$ be color classes of an edge-coloring of $G$. If
\begin{align*}
	\sum_{i \in I} |E_i| \ge r\Delta (|I| - 1) + 1 \hspace{10mm} \text{for every } I \subseteq [n],
\end{align*}
then there exists a full rainbow matching. In particular, if $|E_i| \ge r\Delta$ for every $i$, then there exists a full rainbow matching.
\end{thm}

Our first main result is that the bound $|E_i| \ge r\Delta$ in Theorem \ref{max-degree-thm} is best possible for all $r$ and $\Delta$. 

\begin{thm} \label{hypergraph-construction}
For all integers $r \ge 1$ and $\Delta \ge 2$, there exists an $r$-graph $G$ with maximum degree $\Delta$ and an edge-coloring of $G$ into color classes $E_1, \ldots, E_n$, such that $|E_i| \ge r\Delta - 1$ for every $i$ and there is no full rainbow matching.
\end{thm}

We will give a large family of constructions for Theorem \ref{hypergraph-construction}, with some of our examples being $r$-partite and some not, even for the same values of $r$ and $\Delta$. Our $r$-graphs will be disjoint unions of blow-ups of \textit{nets} from design theory, which are equivalent to sets of mutually orthogonal Latin squares. The number of color classes $n$ will generally be on the order of $r\Delta$. Some of our color classes could have size greater than $r\Delta - 1$, but then it is always possible to obtain constructions where all color classes have size $r\Delta - 1$ by deleting an appropriate selection of edges. For the case of multigraphs $r = 2$, Theorem \ref{hypergraph-construction} was previously known for $\Delta = 2$, but the best constructions for larger values of $\Delta$ had color class sizes $|E_i| = 2 \Delta - 2$ for every $i$. These latter constructions, due to Aharoni, Berger, Kotlar, and Ziv \cite{AhBeKoZi} and to Gao, Ramadurai, Wanless, and Wormald \cite{GaRaWaWo}, disproved earlier conjectures by Aharoni and Berger \cite{AhBe}. Thus, Theorem \ref{hypergraph-construction} improves on previous results for $r = 2$ and $\Delta > 2$ while also generalizing to all larger uniformities $r$. 

The multi-hypergraphs that we construct for Theorem \ref{hypergraph-construction} have many parallel edges. Our next result gives some slightly worse constructions when we assume that the hypergraph edges do not intersect too much. Given an integer $t \ge 1$, a multi-hypergraph is said to be \textit{$t$-simple} if the intersection of any two distinct edges has size at most $t$ (where parallel edges are considered distinct). When $t = 1$, such a multi-hypergraph is a \textit{linear} hypergraph. We prove the following result, which generalizes the case $r = 2$ and $t = 1$ shown by Gao, Ramadurai, Wanless, and Wood \cite{GaRaWaWo}.

\begin{thm} \label{hypergraph-construction-intersecting}
For all integers $1 \le t \le r$ and $\Delta \ge 2$, there exist a $t$-simple, $r$-partite $r$-graph $G$ with maximum degree $\Delta$ and an edge-coloring of $G$ into color classes $E_1, \ldots, E_n$, such that $|E_i| \ge r(\Delta - 1) + t - 1$ for every $i$ and there is no full rainbow matching.
\end{thm}

Next, we focus on full rainbow matchings when we assume that the edge-colorings are proper, i.e. that each color class $E_i$ is a matching. Our constructions for Theorem \ref{hypergraph-construction} and Theorem \ref{hypergraph-construction-intersecting} are far from being proper edge-colorings, as many of the color classes $E_i$ have maximum degree within a constant factor of the maximum degree $\Delta$ of the entire $r$-graph $G$. When we restrict the maximum degree of each color class, optimal conditions for the existence of full rainbow matchings can often be improved. For example, Delcourt and Postle \cite{DePo} the following asymptotic theorem. Recall that the \textit{codegree} of a pair of distinct vertices $u, v$ in a multi-hypergraph is the number of edges containing both $u$ and $v$.

\begin{thm}[Delcourt, Postle] \label{delcourt-postle}
Let $G$ be an $r$-graph with maximum degree $\Delta$, and let $E_1, \ldots, E_n$ be color classes of an edge-coloring of $G$. For fixed $r$ and increasing $\Delta$, if $G$ has maximum codegree $o(\Delta)$, the maximum degree of each $E_i$ is $o(\Delta)$, and $|E_i| \ge (1 + o(1))\Delta$, then there exists a full rainbow matching.
\end{thm}

Delcourt and Postle \cite{DePo} proved more strongly that under the hypotheses of Theorem \ref{delcourt-postle} there are $\Delta$ disjoint full rainbow matchings. Chakraborti and Loh \cite{ChLo} proved a slightly weaker version of Theorem \ref{delcourt-postle} for simple graphs but with an explicit error term. As attempts to optimize their asymptotic result, Delcourt and Postle \cite{DePo} also formulated two conjectures on full rainbow matchings in properly edge-colored multigraphs.

\begin{conj}[Delcourt, Postle] \label{del-pos-conj-1}
Let $G$ be a bipartite multigraph, and let $E_1, \ldots, E_n$ be color classes of a proper edge-coloring of $G$. If $|E_i| \ge \Delta+1$ for every $i$, then there exists a full rainbow matching.
\end{conj}

\begin{conj}[Delcourt, Postle] \label{del-pos-conj-2}
Let $G$ be a general multigraph, and let $E_1, \ldots, E_n$ be color classes of a proper edge-coloring of $G$. If $|E_i| \ge \Delta+2$ for every $i$, then there exists a full rainbow matching.
\end{conj}

These conjectures were suggested as ``sparse" analogues of conjectures of Aharoni and Berger \cite{AhBe} and of Gao, Ramadurai, Wanless, and Wormald \cite{GaRaWaWo} mentioned. This follows the line of various other works that formulate generalizations of well-known conjectures about rainbow matchings and other structures (see, e.g., \cite{AhBe1, AhKo, PoSu}). Here, we disprove both of these conjectures with the following theorem.

\begin{thm} \label{proper-edge-coloring-multigraph-1}
The following statements hold:
\begin{enumerate}
	\item[(1)] For every multigraph $H$ with maximum degree $\Delta \ge 2$ and chromatic index $\chi'$, there exists an associated multigraph $G$ with maximum degree $\Delta$ and chromatic index $\chi'$ and a proper edge-coloring of $G$ into color classes $E_1, \ldots, E_n$, such that $|E_i| \ge \chi' - 1$ for every $i$ and there is no full rainbow matching. If $|E(H)| \le 2\Delta - 1$, then we can take $|E_i| \ge \chi'$ for every $i$.
	\item[(2)] For every $\Delta \ge 2$, there exist a bipartite graph $G$ with maximum degree $\Delta$ and a proper edge-coloring of $G$ into color classes $E_1, \ldots, E_n$, such that $|E_i| \ge \Delta+1$ for every $i$ and there is no full rainbow matching.
	\item[(3)] For every integer $\Delta \ge 3$ with $\Delta \equiv 3, 0 \pmod 4$, there exist a multigraph $G$ with maximum degree and chromatic index $\Delta$, and a proper edge-coloring of $G$ into color classes $E_1, \ldots, E_n$, such that $|E_i| \ge \Delta+2$ for every $i$ and there is no full rainbow matching. If $\Delta \in \{2^m - 1, 2^m\}$ for some integer $m \ge 2$, then we can take $G$ to be a simple graph.
\end{enumerate}
\end{thm}

Applying statement (1) of Theorem \ref{proper-edge-coloring-multigraph-1} with $H$ being the Shannon triangle (the triangle with $\lfloor \Delta/2 \rfloor$ parallel edges on two sides and $\lceil \Delta/2 \rceil$ parallel edges on the third side), we get examples of properly edge-colored multigraphs where $|E_i| \ge \lfloor 3\Delta/2 \rfloor$ for every $i$ and there is no full rainbow matching, which strongly disproves Conjecture \ref{del-pos-conj-2}. Statement (1) of Theorem \ref{proper-edge-coloring-multigraph-1} is proven by relating the chromatic index of $H$ to the non-existence of full rainbow matchings in $G$, and moreover $G$ will include a disjoint union of some number copies of $H$, and thus it will have chromatic index at least $\chi'$. Statements (2) and (3) of Theorem \ref{proper-edge-coloring-multigraph-1} demonstrate that in fact Conjecture \ref{del-pos-conj-1} and Conjecture \ref{del-pos-conj-2} are false even if $G$ has chromatic index $\chi' = \Delta$.

Our final focus is on applying results on full rainbow matchings to list edge-colorings of multi-hypergraphs. This generalizes some of the ideas that go into the proof of statement (1) of Theorem \ref{proper-edge-coloring-multigraph-1}. Given a multi-hypergraph $H$ and a list assignment $L = (L(e) : e \in E(H))$ of colors for $E(H)$, an \textit{$L$-coloring} of $H$ is an edge-coloring $\phi : E(H) \rightarrow \bigcup_{e \in E(H)} L(e)$ such that $\phi(e) \in L(e)$ for every $e \in E(H)$. The \textit{list chromatic index} $\chi_\ell'(H)$ of $H$ is the minimum integer $k \ge 1$ such that for every list assignment $L$ for $E(H)$ with $|L(e)| \ge k$ for all $e \in E(H)$, there exists a \textit{proper} $L$-coloring of $H$. Note that the list chromatic index $\chi_{\ell}'(H)$ is always at least the chromatic index $\chi'(H)$. The celebrated List Edge-Coloring Conjecture (conjectured by various authors, see e.g. \cite{BoHa}) asserts that for every multigraph $H$, we in fact have $\chi_{\ell}'(H) = \chi'(H)$. Famously, Galvin \cite{Ga} proved that the List Edge-Coloring Conjecture holds for all bipartite multigraphs, which in particular proved a conjecture of Dinitz on Latin squares.

\begin{thm}[Galvin] \label{galvin}
For every bipartite multigraph $H$, we have $\chi_{\ell}'(H) = \Delta(H)$.
\end{thm} 

The List Edge-Coloring Conjecture has also been proven to hold for a few other classes of multigraphs such as complete graphs of even degree \cite{HaJa} and prime degree \cite{Sc}. Kahn \cite{Ka2} proved that the List Edge-Coloring Conjecture holds asymptotically for all multigraphs $H$, and Kahn \cite{Ka1} also proved that for every uniform multi-hypergraph $H$ with maximum degree $\Delta$ and maximum codegree $o(\Delta)$, we have $\chi_{\ell}'(H) = (1+o(1))\Delta$.

The problem of finding a proper $L$-coloring of a multi-hypergraph $H$ can easily be translated into the problem of finding a full rainbow matching in an auxiliary edge-colored multi-hypergraph $G$ (see Section \ref{connection-to-rainbow-matchings} for details). This rainbow matchings point of view is useful for studying list edge-colorings in a ``color degree" setting. Given a multi-hypergraph $H$ and list assignment $L$ for $E(H)$, the \textit{maximum color degree} of $(H, L)$ is the maximum, over all colors $c$, of the maximum degree of the edge multiset $E_c = \{e \in E(H) : c \in L(e)\}$. The \textit{maximum color codegree} of $(H, L)$ is defined the same way but with ``maximum degree" replaced by ``maximum codegree". The maximum color degree of the pair $(H, L)$ is in fact the same as the maximum degree of the auxiliary multi-hypergraph $G$, each list $L(e)$ corresponds to a color class on $G$, and the associated edge-coloring of $G$ is a proper edge-coloring. Therefore, we can apply Theorem \ref{max-degree-thm} and Theorem \ref{delcourt-postle} to this multi-hypergraph $G$ and derive the following color degree conditions for the existence of proper $L$-colorings of a multi-hypergraph $H$.

\begin{thm} \label{list-edge-coloring-general}
Let $H$ be an $r$-graph and let $L$ be a list assignment for $E(H)$. If $(H, L)$ has maximum color degree $\Delta$ and $|L(e)| \ge r\Delta$ for every $e \in E(H)$, then there exists a proper $L$-coloring of $H$.
\end{thm}

\begin{thm} \label{list-edge-coloring-asymptotic}
Let $H$ be an $r$-graph and let $L$ be a list assignment for $E(H)$. If $(H, L)$ has maximum color degree $\Delta$, maximum color codegree $o(\Delta)$, and $|L(e)| \ge (1+o(1))\Delta$ for every $e \in E(H)$, then there exists a proper $L$-coloring of $H$.
\end{thm}

Observe that Theorem \ref{list-edge-coloring-asymptotic}, derived from Theorem \ref{delcourt-postle} of Delcourt and Postle \cite{DePo}, is a color degree generalization of Kahn's \cite{Ka1} asymptotic result on the list chromatic index of $r$-graphs. It is then natural to ask whether Galvin's Theorem holds in the color degree setting. Our final result is that this is not the case.

\begin{thm} \label{galvin-color-degree}
For every integer $\Delta \ge 2$, there exist a bipartite graph $H$ and list assignment $L$ for $E(H)$ such that $(H, L)$ has maximum color degree $\Delta$, $|L(e)| = \Delta$ for every $e \in E(H)$, and there is no proper $L$-coloring.
\end{thm}

The above list edge-coloring results were inspired by similar work on list vertex-coloring surrounding a disproved conjecture of Reed \cite{Re} (see remarks in Section \ref{section-proof-of-galvin} for details).

This paper is organized as follows. In Section \ref{construction-section}, we present the general construction method and describe a proof of Theorem \ref{max-degree-thm}. In Section \ref{section-general-edge-colorings}, we prove Theorem \ref{hypergraph-construction} and Theorem \ref{hypergraph-construction-intersecting} on full rainbow matchings in general edge-colored $r$-graphs. In Section \ref{section-proper-edge-colorings}, we prove Theorem \ref{proper-edge-coloring-multigraph-1}, on full rainbow matchings in properly edge-colored multigraphs. In Section \ref{section-list-edge-coloring}, we prove Theorem \ref{galvin-color-degree} on list edge-colorings of bipartite graphs. Finally, in Section \ref{section-questions} we conclude this paper with some questions.

\section{The construction method} \label{construction-section}

In this section, we describe the iterative method we will use to construct edge-colored multi-hypergraphs with no full rainbow matchings. We will then explain how Theorem \ref{max-degree-thm} follows from known results, since a proof has not previously appeared in the literature.

\subsection{Constructing multi-hypergraphs with no full rainbow matchings} \label{construction-subsection}

All of the constructions in this paper, beyond the elementary ones, will be derived by iteratively applying the following simple lemma. It is adapted from \cite{HaWd1}, where the lemma was described in terms of independent transversals.

\begin{lem} \label{join-lemma}
Let $G$ and $H$ be disjoint multi-hypergraphs, and let $\mathcal{P} = \{E_1, \ldots, E_m\}$ and $\mathcal{Q} = \{F_1, \ldots, F_n\}$ be color classes of edge-colorings of $G$ and $H$, respectively. Let $\mathcal{R} = \{E_1', \ldots, E_m',$ $F_1, \ldots, F_{n-1}\}$ be color classes of an edge-coloring of $G \cup H$, where $E_1' \supseteq E_1, \ldots, E_m' \supseteq E_m$ are obtained by distributing each of the edges of $F_n$ into one of $E_1, \ldots, E_m$ arbitrarily. If $G$ has no full rainbow matching with respect to $\mathcal{P}$ and $H$ has no full rainbow matching with respect to $\mathcal{Q}$, then $G \cup H$ has no full rainbow matching with respect to $\mathcal{R}$.
\end{lem}

\begin{proof}
Assume for contradiction that $G \cup H$ has full rainbow matching $\{e_1, \ldots, e_m, f_1, \ldots, f_{n-1}\}$ with respect to $\mathcal{R}$, where $e_i \in E_i'$ for $1 \le i \le m$ and $f_j \in F_j$ for $1 \le j \le n-1$. If $e_i \in E_i$ for every $1 \le i \le m$, then $\{e_1, \ldots, e_m\}$ is a full rainbow matching in $G$ with respect to $\mathcal{P}$, a contradiction. So suppose instead that $e_j \in E_j' \cap F_n$ for some $1 \le j \le m$. Then $\{f_1, \ldots, f_{n-1}, e_j\}$ is a full rainbow matching in $H$ with respect to $\mathcal{Q}$, again a contradiction.
\end{proof}

\begin{figure}
	\centering
	\includegraphics[scale=0.85]{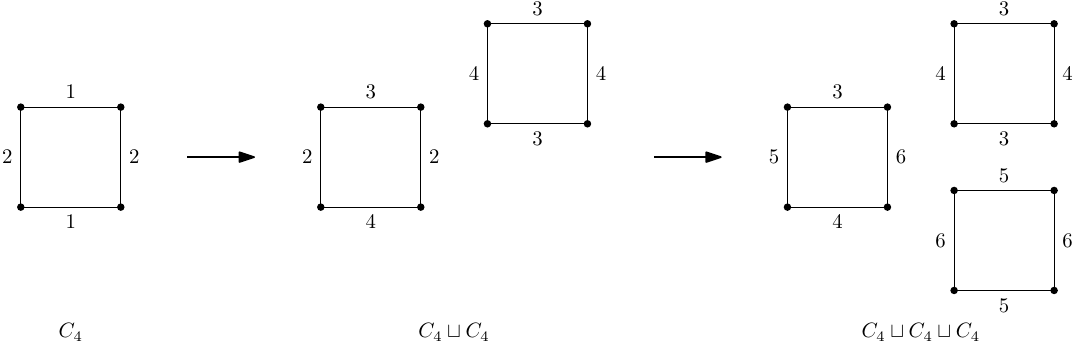}
	\caption{An illustration of iteratively applying Lemma \ref{join-lemma} to the $4$-cycle $C_4$ with a proper $2$-edge-coloring.}
	\label{4-cycles}
\end{figure}

For a basic example of how to apply Lemma \ref{join-lemma}, see Figure \ref{4-cycles}. In this example, we apply Lemma \ref{join-lemma} twice in succession to the $4$-cycle $C_4$ that is properly edge-colored with $2$ colors. In both steps, we choose one color class from the previous step and distribute each of its edges among the newly added color classes. The initial graph $C_4$ has $2$ color classes each of size $2$, and the final edge-colored graph $C_4 \sqcup C_4 \sqcup C_4$ has $4$ color classes each of size $3$. (The symbol $\sqcup$ denotes the disjoint union.) The line graph of this final edge-colored graph is an example of a vertex-partitioned graph with no independent transversal that has been described various times in the literature \cite{HaWd1, Ji, SzTa, Yu}.

Similar to the above example of $C_4$, we will often apply Lemma \ref{join-lemma} repeatedly on a fixed multi-hypergraph $H$ with no full rainbow matching, particularly in Section \ref{section-general-edge-colorings}. We seek to make all the color classes as large as possible, and the following theorem conveniently quantifies how large we can eventually make our color classes if we apply Lemma \ref{join-lemma} repeatedly on a fixed multi-hypergraph.

\begin{thm} \label{join-general-hypergraphs}
Let $H$ be a multi-hypergraph and let $\mathcal{Q} = \{F_1, \ldots, F_m\}$ be color classes of an edge-coloring of $H$, such that $H$ has no full rainbow matching with respect to $\mathcal{Q}$. Suppose that integer $q \ge 1$ satisfies
\begin{align*}
	\sum_{i \in I} |F_i| \ge q(|I| - 1) + 1 \hspace{5mm} \text{ for every } I \subseteq [m].
\end{align*}
Then there exist an integer $k \ge 1$ and an edge-coloring of the multi-hypergraph $H^k = H \sqcup \cdots \sqcup H$ (where $H$ appears $k$ times in the disjoint union) into color classes $\mathcal{P} = \{E_1, \ldots, E_{k(m-1)+1}\}$, such that $|E_i| \ge q$ for every $i$ and there is no full rainbow matching. Specifically, if $J = \{i : |F_i| < q\} \neq \emptyset$ and $F_J = \bigcup_{i \in J} F_i$, then we can take any $k \ge 1 + \sum_{i \in J} \left\lceil \frac{q - |F_i|}{|F_J| - q(|J| - 1)} \right\rceil$.
\end{thm}

See \cite{HaWd1} for a proof of Theorem \ref{join-general-hypergraphs}, which is phrased there in terms of independent transversals but actually works for general simplicial complexes.

\subsection{Proof of Theorem \ref{max-degree-thm}} \label{proof-of-thm-1}
Now we briefly explain how Theorem \ref{max-degree-thm} follows from known results on rainbow matchings. Given a subset $E$ of edges in a multi-hypergraph $G$, the \textit{matching number} $\nu(E)$ is the cardinality of a largest matching in $E$. The \textit{fractional matching number} $\nu^\ast(E)$ is the maximum value of $\sum_{e \in E} x(e)$ over all \textit{fractional matchings} $x : E \rightarrow \mathbb{R}$, which are weight functions that satisfy $x(e) \ge 0$ for all $e \in E(G)$, and $\sum_{e \ni v} x(e) \le 1$ for all $v \in V(G)$. Note that $\nu^\ast(E) \ge \nu(E)$. Aharoni, Berger, and Meshulam \cite{AhBeMe} proved the following Hall-type theorem, which was previously proved by Aharoni and Haxell \cite{AhHa} when $\nu^\ast$ is replaced by $\nu$.

\begin{thm} \label{exists-rainbow-matching}
Let $G$ be an $r$-graph, and let $E_1, \ldots, E_n$ be color classes of an edge-coloring of $G$. If 
\begin{align*}
	\nu^\ast \left( \bigcup_{i \in I} E_i \right) > r(|I| - 1) \hspace{10mm} \text{for all } I \subseteq [n],
\end{align*}
then there exists a full rainbow matching.
\end{thm}

Aharoni, Berger, and Meshulam \cite{AhBeMe} proved Theorem \ref{exists-rainbow-matching} by combining a topological connectedness condition for independent transversals from \cite{AhHa, Me1} together with their own spectral lower bound on the connectedness of the independence complex of a graph. Using Theorem \ref{exists-rainbow-matching}, we deduce Theorem \ref{max-degree-thm} as follows.

\begin{proof}[Proof of Theorem \ref{max-degree-thm}]
Since the $r$-graph $G$ has maximum degree $\Delta$, the constant function $x = 1/\Delta$ is a fractional matching on any edge subset of $E(G)$. Thus, for every $I \subseteq [n]$,
\begin{align*}
	\nu^\ast \left( \bigcup_{i \in I} E_i \right) \ge \frac{1}{\Delta} \left| \bigcup_{i \in I} E_i \right| \ge \frac{1}{\Delta} \left[r\Delta(|I| - 1) + 1\right] > r(|I| - 1).
\end{align*}
By Theorem \ref{exists-rainbow-matching}, there exists a full rainbow matching.
\end{proof}

We note that a non-uniform version of Theorem \ref{max-degree-thm} was stated in \cite{AhBeKoZi} but not proved.

\section{Edge-colored multi-hypergraphs with no full rainbow matchings} \label{section-general-edge-colorings}
In this section, we apply our construction method to derive various edge-colored $r$-graphs with no full rainbow matchings, in particular proving Theorem \ref{hypergraph-construction} and Theorem \ref{hypergraph-construction-intersecting}.

\subsection{Definitions}
Our constructions for Theorem \ref{hypergraph-construction} will be disjoint unions of blow-ups of nets, which are equivalent to sets of mutually orthogonal Latin squares. The main intuition for why nets are useful for us is that their line graphs are balanced complete multipartite graphs. We briefly review nets below as far as they are relevant for our constructions, and see \cite{Ev} for a more extensive overview. 

Fix $r \ge 1$ and $s \ge 2$. An \textit{$s$-net of order $r$}, or an \textit{$(r, s)$-net} for short, is an $r$-graph $H_{r,s}$ with $r^2$ vertices and $rs$ edges, such that the edges $E(H_{r,s})$ can be partitioned into $s$ perfect matchings $F_1, \ldots, F_s$ each of size $r$, and $|e \cap f| = 1$ for all $e \in F_i$ and $f \in F_j$ with $i \neq j$. The matchings $F_1, \ldots, F_s$ are called the \textit{parallel classes} of $H_{r,s}$. For example, an $(r, 2)$-net $H_{r,2}$ is an $r \times r$ grid, with one parallel class $F_1$ being $r$ ``horizontal" edges and the other parallel class $F_2$ being $r$ ``vertical" edges. For $r \ge 2$, the maximum integer $s$ for which there exists some $(r, s)$-net is always at most $r+1$, and if such an $(r, r+1)$-net exists then it is called an \textit{affine plane of order $r$}. An affine plane of order $r$ exists whenever $r$ is a prime power, as it uniquely corresponds to a projective plane of order $r$.

We will only need $(r,s)$-nets for our constructions, but for intuition and understanding when they are $r$-partite, we briefly describe their relationship with mutually orthogonal Latin squares. Two Latin squares $L^1$ and $L^2$ of order $r$ are \textit{orthogonal} if for every pair $(x, y) \in [r] \times [r]$, there is exactly one pair $(i, j) \in [r] \times [r]$ such that $M_{ij}^1 = x$ and $M_{ij}^2 = y$. A set of Latin squares $L^1, \ldots, L^n$ is said to be \textit{mutually orthogonal} if every pair of distinct Latin squares in this set are orthogonal. Then an $(r, s)$-net $H_{r,s}$ is equivalent to a set of $s - 2$ mutually orthogonal Latin squares $L^1, \ldots, L^{s - 2}$. To see this, we start with an $(r, 2)$-net $H_{r,2}$ with parallel classes $F_1 = \{e_1, \ldots, e_r\}$ and $F_2 = \{f_1, \ldots, f_r\}$, and we let $v_{ij}$ denote the unique vertex in $e_i \cap f_j$. Assuming for $2 \le t \le s-1$ that we have constructed a $(r, t)$-net $H_{r, t}$ with parallel classes $F_1, \ldots, F_{t}$, we construct the $r$-graph $H_{r, t+1}$ by adding a parallel class $F_{t+1}$ with edges $\{v_{ij} : L^{t-1}_{ij} = k \}$, for $1 \le k \le r$. The mutual orthogonality of the Latin squares $L^1, \ldots, L^{t-1}$ ensures that $H_{r, t+1}$ is an $(r, t+1)$-net with parallel classes $F_1, \ldots, F_{t+1}$. Moreover, we can reverse this process and derive the $s - 2$ mutually orthogonal Latin squares $L^1, \ldots, L^{s - 2}$ from the $(r, s)$-net $H_{r,s}$. 

For the Latin squares $L^1, \ldots, L^{s-2}$ above, notice that for $2 \le t \le s - 1$, the Latin square $L^{t - 1}$ corresponds to an $r$-partition $V_1, \ldots, V_r$ of the vertices of the $r$-graph $H_{r,t}$ that shows that $H_{r,t}$ an $r$-partite $r$-graph, namely $V_k = \{ v_{ij} : L_{ij}^{t-1} = k \}$ for $1 \le k \le r$. Conversely, any proper $r$-partition $V_1, \ldots, V_r$ of the vertices of $H_{r,t}$ is a decomposition of each of the Latin squares $L^1, \ldots, L^{t-2}$ into transversals, and this is equivalent to a Latin square orthogonal to each of $L^1, \ldots, L^{t-2}$. Therefore, the $(r,s)$-net $H_{r,s}$ is not $r$-partite if and only if the corresponding set of mutually orthogonal Latin squares $L^1, \ldots, L^{s-2}$ is maximal. In this case, we will say that $H_{r,s}$ is a \textit{maximal} net (so maximality is the same as not being $r$-partite).

\subsection{Proof of Theorem \ref{hypergraph-construction}}

Now we construct examples for Theorem \ref{hypergraph-construction}. Example \ref{example-1} below gives our most basic construction that works for all integers $r \ge 1$ and $\Delta \ge 2$, but we will present a larger family of examples together with it.

Fix integers $r \ge 1$ and $s \ge 2$, and suppose that there exists an $(r,s)$-net $H_{r,s}$ with parallel classes $F_1, \ldots, F_s$. (Note that an $(r,s)$-net exists for all $r \ge 1$ and $s \in \{2, 3\}$.) Naturally, we view $F_1, \ldots, F_s$ as the color classes of a proper edge-coloring of $H_{r,s}$. Then every rainbow matching in $H_{r,s}$ has size at most one, since any two edges from distinct color classes intersect in exactly one vertex. Given a sequence of integers $a_1, \ldots, a_s \ge 1$, let $H_{r,s}(a_1, \ldots, a_s)$ denote the $r$-graph obtained from $H_{r,s}$ by replacing, for $1 \le i \le s$, every edge in the parallel class $F_i$ by $a_i$ parallel edges. Let $G_{r,s}(a_1, \ldots, a_s)$ be the disjoint union of $s - 1$ copies of $H_{r,s}(a_1, \ldots, a_s)$, and for $1 \le i \le s$, let $F_i'$ be the set edges of $G_{r,s}(a_1, \ldots, a_s)$ that are copies of edges from $F_i$. Then $F_1', \ldots, F_s'$ are color classes of an edge-coloring of $G_{r,s}(a_1, \ldots, a_s)$, and moreover there is no full rainbow matching, since there are $s$ color classes on $s - 1$ components, and in a given component any two edges from distinct color classes intersect in exactly one vertex.

Note that the $r$-graph $G_{r,s}(a_1, \ldots, a_s)$ has maximum degree $\sum_{i=1}^s a_i$, and that the color class sizes are $|F_i'| = r(s-1)a_i$ for $1 \le i \le s$. We wish to apply Theorem \ref{join-general-hypergraphs} to $G_{r,s}(a_1, \ldots, a_s)$ to obtain an $r$-graph with maximum degree $\Delta$ and color classes $E_1, \ldots, E_n$, such that $|E_i| \ge r\Delta - 1$ for every $i$ and there is no full rainbow matching. For this to work, we require that the sequence of integers $a_1, \ldots, a_{s} \ge 1$ satisfy the following conditions: 
\begin{itemize}
	\item[(1)] the equality $\sum_{i = 1}^{s} a_i = \Delta$, and
	\item[(2)] the inequalities $\sum_{i \in I} r(s-1)a_i > (r\Delta - 1)(|I| - 1)$ for all $I \subseteq [s]$.
\end{itemize}
We present various examples that satisfy these conditions, leaving out the algebraic verifications. We let $G_{r,s}$ denote the final $r$-graph obtained in each case.

\begin{figure}
	\centering
	\includegraphics[scale=0.85]{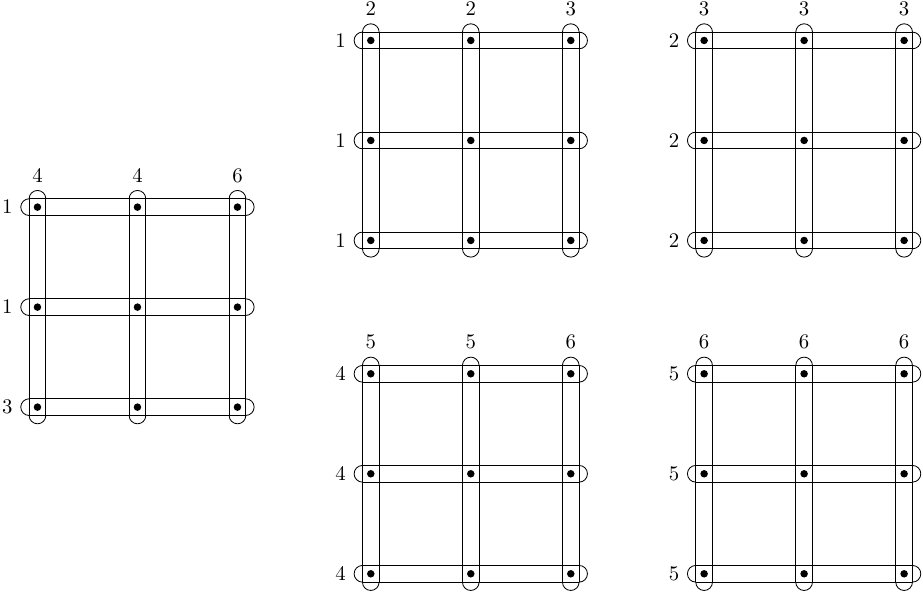}
	\caption{The disjoint union of five copies of the $(3,2)$-net, i.e. the $3 \times 3$ grid, edge-colored into color classes each of size $3 \cdot 2 - 1 = 5$ with no full rainbow matching.}
	\label{3-3-grids}
\end{figure}

\begin{example} \label{example-1}
One example that works for all $r \ge 1$ and $\Delta \ge 2$ (and thus proves Theorem \ref{hypergraph-construction}) is taking $s = 2$ and any two integers $a_1, a_2 \ge 1$ such that $a_1 + a_2 = \Delta$. In this case, Theorem \ref{join-general-hypergraphs} implies that if $G_{r,2}$ is the disjoint union of $k = r\Delta - 1$ copies of $G_{r,2}(a_1, a_2) = H_{r,2}(a_1, a_2)$, then there exists an edge-coloring of $G_{r,2}$ into $m = r \Delta$ color classes $E_1, \ldots, E_{r\Delta}$ such that $|E_i| \ge r\Delta - 1$ for every $i$ and there is no full rainbow matching. This example is always $r$-partite because there always exists a Latin square of order $r$. Note that for $r = \Delta = 2$, the graph $G_{2,2}$ is the disjoint union three $4$-cycles presented from Figure \ref{4-cycles}. For $r = 3$ and $\Delta = 2$, the $3$-graph $G_{3,2}$ and an appropriate edge-coloring is shown in Figure \ref{3-3-grids}.
\end{example}

\begin{example} \label{example-2}
For $r \ge 1$ and $s \ge 2$, suppose that an $(r, s)$-net $H_{r,s}$ exists. Then whenever $\Delta$ is a positive multiple of $s$, we can take $a_1 = \cdots = a_s = \frac{\Delta}{s}$. In this case, Theorem \ref{join-general-hypergraphs} states that if $G_{r,s}$ is the disjoint union of $k = 1 + s \left\lceil \frac{r\Delta - s}{s(s - 1)} \right\rceil$ copies of $G_{r,s}(\frac{\Delta}{s}, \ldots, \frac{\Delta}{s})$, then there exists an edge-coloring of $G_{r,s}$ into $n = s + s(s - 1) \left\lceil \frac{r\Delta - s}{s(s - 1)} \right\rceil$ color classes $E_1, \ldots, E_n$, such that $|E_i| \ge r\Delta - 1$ for every $i$ and there is no full rainbow matching. This example $G_{r,s}$ is $r$-partite if and only if the net $H_{r,s}$ is not maximal.
\end{example}

\begin{example} \label{example-3}
For $r \ge 1$ and $s \ge 2$, assume that an $(r, s)$-net $H_{r,s}$ exists. The following two examples again are either $r$-partite or not $r$-partite depending on whether the net $H_{r,s}$ is maximal.

The first example is that if $\Delta \equiv \ell \pmod{s}$, $1 \le \ell \le s - 1$, and $\Delta > (s - 1)(s - \ell) - \frac{s(s - 1 - \ell)}{r\ell}$, then we can take $a_1 = \cdots = a_{\ell} = \left\lceil \frac{\Delta}{s} \right\rceil$ and $a_{\ell + 1} = \cdots = a_s = \left\lfloor \frac{\Delta}{s} \right\rfloor$. If $\ell = 0$ and $\Delta \ge s$, then this sequence $a_1, \ldots, a_s$ is the same as in Example \ref{example-2}. 

The second example is that if $\Delta - 1 \equiv \ell \pmod{s - 1}$, $0 \le \ell \le s - 2$, and $\Delta > (\ell + 1)(s - 1)$, then we can take $a_1 = \cdots = a_{s-1} = \left\lfloor \frac{\Delta - 1}{s - 1} \right\rfloor$ and $a_s = \Delta - (s - 1) \left\lfloor \frac{\Delta - 1}{s - 1} \right\rfloor$. These are the examples we wanted to show for Theorem \ref{hypergraph-construction}.
\end{example}

There are many other integer sequences $a_1, \ldots, a_s \ge 1$ that satisfy conditions (1) and (2), and many more constructions are possible using nets. It is then natural to ask whether there are constructions for Theorem \ref{hypergraph-construction} that do not involve blow-ups of nets as components. The answer is yes when $r = \Delta = 2$. Recall that in this case the construction from Example \ref{example-1} is the disjoint union of three $4$-cycles, which was illustrated in Figure \ref{4-cycles}. It was shown in \cite{HaWd1} that for any disjoint union of three cycles each of length $1$ modulo $3$, there exists an edge-coloring of it into color classes $E_1, \ldots, E_n$ such that $|E_i| = 3$ for every $i$ and there is no full rainbow matching. These disjoint unions are the unique minimal such multigraphs with maximum degree $2$, by a result of Aharoni, Holzman, Howard, and Spr\"ussel \cite{AhHoHoSp}. It would be interesting to find other values of $r$ and $\Delta$ for which there exist constructions for Theorem \ref{hypergraph-construction} not involving blow-ups of nets.

Finally, we observe that our large family of extremal constructions for Theorem \ref{hypergraph-construction} stands in contrast to extremal constructions for independent transversals in graphs. While a great variety of multi-hypergraphs can be used to achieve the color class sizes $r\Delta - 1$ with no full rainbow matching as we have seen, the kinds of graphs that achieve the analogous class sizes $2\Delta - 1$ with no independent transversal (this bound is best possible by a theorem of Haxell \cite{Ha2}) are always disjoint unions of the complete bipartite graph $K_{\Delta,\Delta}$ (see \cite{AhHoHoSp, HaWd1, Ji, SzTa, Yu}). It is also interesting that both $r$-partite and non-$r$-partite $r$-graphs can be used to achieve color class sizes $r\Delta - 1$ with no full rainbow matching.

\subsection{When the number of color classes is small}
If the number of color classes $n$ is slightly less than $r\Delta$, then the sufficient condition $|E_i| \ge r\Delta$ in Theorem \ref{max-degree-thm} can be improved. Specifically, the more general condition stated in Theorem \ref{max-degree-thm} implies the following.

\begin{thm} \label{max-degree-number-of-classes}
Let $G$ be an $r$-graph with maximum degree $\Delta$, and let $E_1, \ldots, E_n$ be color classes of an edge-coloring of $G$, where $n \ge 2$. If $|E_i| \ge \left\lfloor \frac{(n - 1)r\Delta}{n} \right\rfloor + 1$ for every $i$ then there exists a full rainbow matching.
\end{thm} 

Here, we show that the bound $|E_i| \ge \left\lfloor \frac{(n - 1)r\Delta}{n} \right\rfloor + 1$ in Theorem \ref{max-degree-number-of-classes} is best possible for certain values of $r$, $\Delta$, and $n$. 

\begin{thm}
Fix integers $r \ge 1$ and $s \ge 2$, suppose that an $(r, s)$-net exists, $\Delta \equiv 0 \pmod{s}$, and $n \equiv s \pmod{s(s - 1)}$. There exist $r$-graphs $G$ with maximum degree $\Delta$ and edge-colorings into color classes $E_1, \ldots, E_n$ such that $|E_i| \ge \left\lfloor \frac{(n - 1)r\Delta}{n} \right\rfloor$ for every $i$ and there is no full rainbow matching.
\end{thm}

\begin{proof}
Like in Example \ref{example-2}, we apply Theorem \ref{join-general-hypergraphs} to the $r$-graph $G_{r,s}(\frac{\Delta}{s}, \ldots, \frac{\Delta}{s})$ and the color classes $F_1', \ldots, F_s'$ each of size $|F_i'| = \frac{r(s-1)\Delta}{s}$, but this time when we apply Theorem \ref{join-general-hypergraphs}, we take $q = \left\lfloor \frac{(n - 1)r\Delta}{n} \right\rfloor$. We can write $n \ge \frac{r\Delta}{r\Delta - q} = n_0$. Theorem \ref{join-general-hypergraphs} says that for some $k \ge 1$ we can find an edge-coloring of the disjoint union of $k$ copies of $G_{r,s}(\frac{\Delta}{s}, \ldots, \frac{\Delta}{s})$ with $k(s-1)+1$ color classes $E_1, \ldots, E_{k(s-1)+1}$ each of size at least $q$, such that there is no full rainbow matching. Specifically, we can take any integer $k \ge 1+s\left\lceil \frac{q - \frac{r(s-1)\Delta}{s}}{r(s-1)\Delta - q(s-1)} \right\rceil = 1+s\left\lceil \frac{qs - r(s-1)\Delta}{s(s-1)(r\Delta - q)} \right\rceil = 1+s\left\lceil \frac{n_0 - s}{s(s-1)} \right\rceil$. In particular, we can take $k = 1+s\left\lceil \frac{n - s}{s(s-1)} \right\rceil = 1+\frac{n-s}{s-1}=\frac{n-1}{s-1}$, and then the number of color classes is $k(s-1)+1=n$, as required.
\end{proof}

\subsection{Proof of Theorem \ref{hypergraph-construction-intersecting}}
Now we give our construction for Theorem \ref{hypergraph-construction-intersecting}. It is a modification of the $r$-graph $G_{r,2}$ from Example \ref{example-1}, which was a disjoint union of blow-ups of the $r \times r$ grid $H_{r,2}$. 

Fix $1 \le t \le r$. For $d \ge 1$, the \textit{$(r,t,d)$-sunflower} $S_{r,t,d}$ is the $r$-graph with $d(r - t) + t$ vertices and $d$ edges, such that there exists a vertex subset $K \subseteq V(S_{r,t,d})$ of size $t$ such that $e \cap f = K$ for any two distinct edges $e, f \in E(S_{r,t,d})$. The set $K$ is called the \textit{kernel} of the sunflower $S_{r,t,d}$. Note that $S_{r,t,d}$ is a $t$-simple, $r$-partite $r$-graph with maximum degree $d$. For $\Delta \ge 2$, let $H_{r,t,\Delta}$ be the $t$-simple, $r$-partite $r$-graph with maximum degree $\Delta$ constructed as follows. First we take the disjoint union of $r$ copies $S_1, \ldots, S_r$ of the sunflower $S_{r,t,\Delta-1}$, and say $S_i$ has kernel $K_i = \{v_{i1}, \ldots, v_{it}\}$ for $1 \le i \le r$. Then we add the edges $e_j = \{v_{1j}, \ldots, v_{rj}\}$ for $1 \le j \le t$, which form a perfect matching on $K_1 \cup \cdots \cup K_r$. Let $F_1 = \bigcup_{i \in [r]} E(S_i)$ and let $F_2 = \{e_1, \ldots, e_t\}$. Then $F_1, F_2$ form the color classes of an edge-coloring of $H_{r,t,\Delta}$ with no full rainbow matching.

\begin{example} \label{example-5} 
We apply Theorem \ref{join-general-hypergraphs} to the multi-hypergraph $H_{r,t,\Delta}$ with the color classes $F_1, F_2$. Note that $|F_1| = r(\Delta - 1)$ and $|F_2| = t$. Then applying Theorem \ref{join-general-hypergraphs} with $q = r(\Delta - 1) + t - 1$, if $G_{r,t,\Delta}$ is the disjoint union of $k = r(\Delta - 1)+t-1$ copies of $H_{r,t,\Delta}$, there exists an edge-coloring of $G_{r,t,\Delta}$ into $k+1 = r(\Delta - 1)+t$ color classes $E_1, \ldots, E_{r(\Delta-1)+t}$ such that $|E_i| \ge r(\Delta - 1) + t - 1$ for every $i$ and there is no full rainbow matching. Since $G_{r,t,\Delta}$ is a $t$-simple, $r$-partite $r$-graph with maximum degree $\Delta$, this completes the proof of Theorem \ref{hypergraph-construction-intersecting}. 
\end{example}

\section{Properly edge-colored multigraphs with no full rainbow matchings} \label{section-proper-edge-colorings}

In this section, we apply our construction method to derive properly edge-colored multigraphs with no full rainbow matchings, which prove Theorem \ref{proper-edge-coloring-multigraph-1}. Before proving these theorems, we outline some (non-optimal) examples of properly edge-colored multigraphs with no full rainbow matchings that we will use.

\subsection{Preliminary examples} \label{preliminary-examples}
For proving statement (1) of Theorem \ref{proper-edge-coloring-multigraph-1}, we will use the following basic construction.

\begin{prop} \label{chromatic-index-rainbow-matching}
Let $H$ be a multigraph with chromatic index $\chi'$, and let $G$ be the disjoint union of $\chi' - 1$ copies of $H$. Assuming that the edge set of $H$ is $\{e_1, \ldots, e_m\}$, consider the proper edge-coloring of $G$ into color classes $E_1, \ldots, E_m$ where $E_i$ consists of the copies of the edge $e_i$ among the disjoint copies of $H$ in $G$. Then there is no full rainbow matching.
\end{prop}

\begin{proof}
A full rainbow matching in $G$ corresponds to a proper edge-coloring of $H$ using $\chi' - 1$ colors: If $H_1, \ldots, H_{\chi' - 1}$ are the copies of $H$ in $G$, then a full rainbow matching $M$ of $G$ corresponds to the proper edge-coloring of $H$ where the edges of $H$ whose copies lie in $M \cap E(H_j)$ are given the color $j$, for $j \in \{1, \ldots, \chi' - 1\}$. Since there is no proper edge-coloring of $H$ using $\chi' - 1$ colors, there is no full rainbow matching in $G$.
\end{proof}

Next, for proving statements (2) and (3) of Theorem \ref{proper-edge-coloring-multigraph-1}, we will utilize some known constructions from the literature. Recall that Aharoni and Berger \cite{AhBe} conjectured the following.

\begin{conj}[Aharoni, Berger] \label{aharoni-berger-bipartite}
Let $G$ be a bipartite multigraph, and let $E_1, \ldots, E_n$ be color classes of a proper edge-coloring of $G$. If $|E_i| \ge n+1$ for every $i$, then there exists a full rainbow matching.
\end{conj}

The bound $|E_i| \ge n+1$ in Conjecture \ref{aharoni-berger-bipartite} cannot be improved for even $n \ge 2$, as shown by the following well-known construction which we reprove for convenience.

\begin{prop} \label{latin-square-no-transversal}
For every even integer $n \ge 2$, there exists a proper edge-coloring of $G = K_{n,n}$ into $n$ perfect matchings each of size $n$, such that there is no full rainbow matching.
\end{prop}

\begin{proof}
Label the vertices of each part of $K_{n,n}$ as $A = \{x_0, \ldots, x_{n-1}\}$ and $B = \{y_0, \ldots, y_{n-1}\}$. Taking the additive abelian group $\mathbb{Z}_n$ to be our color set, consider the edge-coloring of $K_{n,n}$ where edge $x_iy_j$ is given color $i+j$, for all $i, j \in \mathbb{Z}_n$. It is easy to see that this edge-coloring is proper. Suppose for contradiction that there exists a full rainbow matching $M$, where edge $x_{a_k}y_{b_k} \in M$ is given color $k \in \mathbb{Z}_n$. Then modulo $n$,
\begin{align*}
	\sum_{k \in \mathbb{Z}_n} a_k = \sum_{k \in \mathbb{Z}_n} b_k = \sum_{k \in \mathbb{Z}_n} k = \frac{n(n - 1)}{2}
\end{align*}
because the matching is perfect. Thus, $\sum_{k \in \mathbb{Z}_n} (a_k + b_k) = n(n-1) = 0$. But also $\sum_{k \in \mathbb{Z}_n} (a_k + b_k) = \sum_{k \in \mathbb{Z}_n} k = \frac{n(n - 1)}{2}$ because the matching is full rainbow, and this is not $0$ when $n$ is even, a contradiction.
\end{proof}

Proposition \ref{latin-square-no-transversal} equivalently states that for all even $n \ge 2$, there exists a Latin square of order $n$ with no transversal. The proof above gives the example of the Cayley table of the additive group $\mathbb{Z}_n$. The Ryser-Brualdi-Stein Conjecture asserts that for odd $n$, there cannot be a proper edge-coloring of $K_{n,n}$ into $n$ perfect matchings that has no full rainbow matching. We note, however, that one can still show that the condition $|E_i| \ge n+1$ in Conjecture \ref{aharoni-berger-bipartite} would be best possible also for odd integers $n \ge 3$, by slightly modifying a multigraph construction of Bar\'at and Wanless \cite{BaWa} (who extended a construction of Drisko \cite{Dr}).

Next, recall that Gao, Ramadurai, Wanless, and Wormald \cite{GaRaWaWo} conjectured the following.

\begin{conj}[Gao, Ramadurai, Wanless, Wormald] \label{aharoni-berger-nonbipartite}
Let $G$ be a multigraph, and let $E_1, \ldots, E_n$ be color classes of a proper edge-coloring of $G$. If $|E_i| \ge n + 2$ for every $i$, then there exists a full rainbow matching.
\end{conj}

The bound $|E_i| \ge n+2$ in Conjecture \ref{aharoni-berger-nonbipartite} also cannot be improved for many values of $n$. For example, if $n = 3$ then take $G$ to be the disjoint union of two copies of the complete graph $K_4$ that is properly edge-colored into $3$ perfect matchings each of size $4$, which has no full rainbow matching. (This is the construction of Example \ref{example-2} when $r = 2$ and $s = \Delta = 3$.) One way to generalize this example to other values of $n$ is to use, in place of $K_4$, a construction of Bar\'at, Gy\'arf\'as, and S\'ark\"ozy \cite{BaGySa} (see also \cite{AhBeChHoSe}), as follows.

\begin{prop} \label{barat-gyarfas-sarkozy}
For every integer $n \ge 3$ with $n \equiv 3 \pmod{4}$, there exist a multigraph $G$ and a proper edge-coloring of $G$ into $n$ perfect matchings each of size $n+1$, such that there is no full rainbow matching.
\end{prop}

\begin{proof}
Let $n = 2m-1$. Let $C$ be the cycle of length $2m$ where every edge is replaced with $m-1$ parallel edges. We start with a proper edge-coloring of $C$ into $2m - 2$ perfect matchings, each of size $m$. Let $v_1, \ldots, v_{2m}$ be the vertices of $C$ in cyclic order. Let $H$ be obtained from $C$ by adding a perfect matching $M$ all of whose edges $v_iv_j$ have $j - i$ even (e.g. $M = \{v_{4k-3}v_{4k-1}, v_{4k-2}v_{4k} : 1 \le k \le \frac{m}{2}\}$). We give the edges of $M$ a new color, so that $H$ is properly edge-colored into $2m-1$ perfect matchings each of size $m$. We claim that every rainbow matching in $H$ has size at most $m-1$. This claim follows from the observations that every perfect matching of $C$ has two edges of the same color, and that for any edge $e \in M$ there is no perfect matching of $C \cup \{e\}$ containing $e$. Now let $G$ be the disjoint union of two copies of $H$ each with this proper edge-coloring. Then $G$ is properly edge-colored into $n=2m-1$ color classes each of size $n+1=2m$, and there is no full rainbow matching.
\end{proof}

The above example of Bar\'at, Gy\'arf\'as, and S\'ark\"ozy \cite{BaGySa} has many parallel edges. Now we present a different generalization of the above example of the disjoint union of two copies of $K_4$. This construction has not previously appeared in the literature as far as we are aware.

\begin{prop} \label{complete-graph-rainbow}
For every integer $m \ge 2$, there exists a proper edge-coloring of the disjoint union $G$ of two copies of the complete graph $K_{2^m}$ into $2^m - 1$ perfect matchings each of size $2^m$, such that there is no full rainbow matching.
\end{prop}

\begin{proof}
For the complete graph $K_{2^m}$, take its vertex set to be the additive elementary abelian group $\mathbb{Z}_2^m$. Take $\mathbb{Z}_2^m \setminus \{(0,\ldots,0)\}$ to be our color set, and consider the edge-coloring of $K_{2^m}$ where for two vertices $x, y \in V(K_{2^m})$ we give the edge $xy$ the color $x+y$. It is easy to see that this edge-coloring is proper. Letting $G = K_{2^m} \sqcup K_{2^m}$, assign this proper edge-coloring to both copies of $K_{2^m}$ in $G$. Now suppose for contradiction that there exists a full rainbow matching $M = \{x_1y_1, \ldots, x_{2\ell-1} y_{2\ell-1}\}$ in $G$, where $\ell = 2^{m-1}$. (We slightly abuse notation and consider $\mathbb{Z}_2^m$ as the vertex set of both components of $G$). Since the maximum size of a matching in either component of $G$ is $\ell = 2^{m-1}$, there is one subset $M_1 = \{x_1y_1, \ldots, x_{\ell}y_{\ell}\} \subset M$ of size $\ell$ whose edges lie in one component of $G$, and there is another subset $M_2 = M \setminus M_1 = \{x_{\ell+1}y_{\ell+1}, \ldots, x_{2\ell - 1}y_{2\ell - 1}\}$ of size $\ell - 1$ whose edges lie in the other component of $G$. The matching $M_1$ is a perfect matching on its component, so $\sum_{i=1}^\ell (x_i + y_i) = 2^{m-1}(1, \ldots, 1) = (0, \ldots, 0)$. On the other hand, the matching $M_2$ covers all but two vertices on its component, and because these two uncovered vertices have a nonzero sum, we have $\sum_{i=\ell+1}^{2\ell-1} (x_i + y_i) \neq (0, \ldots, 0)$. It follows that $\sum_{i=1}^{2\ell-1} (x_i + y_i) \neq (0, \ldots, 0)$. But we also have that $\sum_{i=1}^{2\ell-1} (x_i + y_i) = (0, \ldots, 0)$ because the matching is full rainbow, a contradiction.
\end{proof}

\subsection{Applying the construction method}
Now we prove each statement of Theorem \ref{proper-edge-coloring-multigraph-1}. Notice that examples in Section \ref{section-proper-edge-colorings} already come quite close to proving each of the statements. The goal now is to increase the size of each color class by one or two while maintaining the properties of being a proper edge-coloring and having no full rainbow matching. In some cases, it suffices to apply Theorem \ref{join-general-hypergraphs} to these properly edge-colored multigraphs, similar to what was done in Section \ref{section-general-edge-colorings}. However, we will use a different strategy that works more often in the proper edge-coloring setting. Roughly, we will start with an improperly edge-colored multigraph with large color classes, and then we will derive a properly edge-colored multigraph from it (with necessarily smaller color classes) by taking the disjoint union of it with one of our examples above and applying Lemma \ref{join-lemma}.


\begin{proof}[Proof of Theorem \ref{proper-edge-coloring-multigraph-1}, Statement (1)]
Recall from Example \ref{example-1} that the multigraph $G_{2,2}$ is the disjoint union of $2\Delta - 1$ subgraphs $H_1, \ldots, H_{2\Delta-1}$ each of which is a  blow-up of a $4$-cycle, specifically a $4$-cycle where in one of the pairs of non-adjacent edges we replace each edge by $\Delta-1$ parallel edges. We saw that there is an edge-coloring of $G_{2,2}(1, \Delta-1)$ into $2\Delta$ color classes $F_1, \ldots, F_{2\Delta}$ such that $|F_i| = 2\Delta - 1$ for every $i$ and there is no full rainbow matching. Note that $G_{2,2}(1, \Delta-1)$ is a bipartite multigraph with maximum degree and chromatic index $\Delta$.
	
Now, let $H$ be a multigraph from the hypothesis of the theorem statement with $2 \le \Delta \le \chi' \le m = |E(H)|$. Then let $G'$ be the multigraph from Proposition \ref{chromatic-index-rainbow-matching} with color classes $\mathcal{P}' = \{E_1', \ldots, E_m'\}$, where $|E_i'| = \chi' - 1$ for every $i$. If $m \ge 2\Delta$, then $G', \mathcal{P}'$ is already the desired properly edge-colored multigraph. Assume instead that $m \le 2\Delta - 1$. Then we start with the multigraph $G_{2,2}(1, \Delta - 1)$ and color classes $F_1, \ldots, F_{2\Delta}$. Iteratively for each color class $F_i$, $i \in \{1, \ldots, 2\Delta\}$, we apply Lemma \ref{join-lemma} by adding a disjoint copy of $G', \mathcal{P}'$ and distributing each of the edges of $F_i$ into one of the classes of $\mathcal{P}'$. Specifically, since $F_i$ is a disjoint union of collections of parallel edges with maximum degree at most $\Delta$ and $|F_i| = 2\Delta-1 \ge m$, we partition $F_i$ into $m$ nonempty matchings, and we distribute each of these matchings into a unique color class $E_i'$ in $\mathcal{P}'$. After doing this for every color class $F_i$, the result is a multigraph $G$ that is properly edge-colored into $n = 2\Delta m$ color classes $E_1, \ldots, E_n$, with $|E_i| \ge \chi'$ for every $i$ and there is no full rainbow matching.
\end{proof}

\begin{proof}[Proof of Theorem \ref{proper-edge-coloring-multigraph-1}, Statement (2)]
Recall from Example \ref{example-5} that the graph $G_{2,1,\Delta}$ is the disjoint union of $2\Delta - 2$ subgraphs $H_1, \ldots, H_{2\Delta - 2}$ each of which is a double star, that is, each $H_i$ is a tree with two central adjacent vertices $v_i$ and $w_i$, and $v_i$ also connected to $\Delta - 1$ leaves with edges $e_{i,1}, \ldots, e_{i,\Delta-1}$, and $w_i$ also connected to $\Delta - 1$ leaves with edges $f_{i,1}, \ldots, f_{i,\Delta-1}$. The color classes of $G_{2,1,\Delta}$ are $F_1, \ldots, F_{2\Delta-1}$, where $F_{2\Delta-1} = \{v_1w_1, \ldots, v_{2\Delta-2}w_{2\Delta-2}\}$ and $F_i = \{e_{i,1}, \ldots, e_{i,\Delta-1}, f_{i,1}, \ldots, f_{i, \Delta-1}\}$ for $i \in \{1, \ldots, 2\Delta-2\}$. Then  $|F_i| = 2\Delta-2$ for every $i$, and there is no full rainbow matching. Note that $G_{2,1,\Delta}$ is a bipartite simple graph with maximum degree and chromatic index $\Delta$.

We start with the graph $G_{2,1,\Delta}$ and color classes $F_1, \ldots, F_{2\Delta-1}$. If $\Delta$ is even, then from Proposition \ref{latin-square-no-transversal} let $\mathcal{P}' = \{E_1', \ldots, E_\Delta'\}$ be the color classes of a proper edge-coloring of the complete bipartite graph $G' = K_{\Delta,\Delta}$ where $|E_i'| = \Delta$ for every $i$ and there is no full rainbow matching. Iteratively for each of the color classes $F_i$, $i \in \{1, \ldots, 2\Delta-1\}$, we apply Lemma \ref{join-lemma} by adding a disjoint copy of $G', \mathcal{P}'$ and distributing each of the edges of $F_i$ into one of the classes of $\mathcal{P}'$. Specifically, since each $F_i$ for $i \in \{1,\ldots,2\Delta-2\}$ is the disjoint union of two copies of the star $K_{1,\Delta-1}$, while $F_{2\Delta-1}$ is already a matching, we partition each $F_i$ into $\Delta$ nonempty matchings, and we distribute each of these matchings into a unique color class $E_i'$ in $\mathcal{P}'$. After doing this for every such class $F_i$, the result is a bipartite simple graph $G$ that is properly edge-colored into $n = (2\Delta-1)\Delta$ color classes $E_1, \ldots, E_n$, such that $|E_i| \ge \Delta+1$ for every $i$ and there is no full rainbow matching.

If $\Delta$ is odd, we instead use the graph $G' = K_{\Delta-1,\Delta-1}$ and color classes $\mathcal{P}' = \{E_1', \ldots, E_{\Delta-1}'\}$ from Proposition \ref{latin-square-no-transversal}, with $|E_i'| = \Delta-1$ for every $i$. For each of the classes $F_i$, $i \in \{1, \ldots, 2\Delta-1\}$, of $G_{2,1,\Delta}$ we partition it into $\Delta-1$ matchings each of size $2$, and we distribute the edges of each of these matchings into a unique color class $E_i' \in \mathcal{P}'$. Then the sizes of the color classes $E_i'$ increase from $\Delta - 1$ to $\Delta+1$.
\end{proof}

\begin{proof}[Proof of Theorem \ref{proper-edge-coloring-multigraph-1}, Statement (3)]
The proof follows the same structure as the proof of Theorem \ref{proper-edge-coloring-multigraph-1}, Statement (2), but in place of Proposition \ref{latin-square-no-transversal}, we use Proposition \ref{barat-gyarfas-sarkozy} for the case of multigraphs, and Proposition \ref{complete-graph-rainbow} for the case of simple graphs.
\end{proof}

\section{List edge-coloring and maximum color degree} \label{section-list-edge-coloring}
In this section, we will apply our construction method to list edge-coloring. Specifically, we will prove Theorem \ref{galvin-color-degree} which demonstrates that a color degree generalization of Galvin's theorem \cite{Ga} does not hold.

\subsection{Connection between list edge-colorings and rainbow matchings} \label{connection-to-rainbow-matchings}
Let $H$ be a multi-hypergraph and let $L = (L(e) : e \in E(H))$ be a list assignment for $E(H)$. We explain how to transform the problem of finding a proper $L$-coloring of a multi-hypergraph $H$ into the problem of finding a full rainbow matching in an auxiliary properly edge-colored multi-hypergraph $G$. This is an adaptation of a similar argument for list vertex-coloring \cite{Ha2}. 

For a color $c \in \bigcup_{e \in E(H)} L(e)$, we let $H_c$ be the sub-hypergraph of $H$ with edge set $E_c = \{e \in E(H) : c \in L(e)\}$ and vertex set $V_c = \bigcup_{e \in E_c} e$ (so that $H_c$ has no isolated vertices, for convenience). In $H_c$ we relabel the copy of edge $e \in E(H)$ by $(e, c)$, and we let $G$ be the disjoint union of the sub-hypergraphs $H_c$ over all colors $c$. Assuming that $H$ has edge set $E(H) = \{e_1, \ldots, e_n\}$, consider the edge-coloring of $G$ with color classes $E_1, \ldots, E_n$, where $E_i = \{ (e_i, c) : c \in L(e_i) \}$ for all $i$. Note that each color class $E_i$ is a matching since its edges lie in different sub-hypergraphs $H_c$. That is, this edge-coloring of $G$ is proper. Then a proper $L$-coloring $\phi$ of $H$ corresponds to a unique full rainbow matching in $G$, namely $\{(e, \phi(e)) : e \in E(H)\}$. Recall that the \textit{maximum color degree} of $(H, L)$ is the maximum, over all colors $c$, of the maximum degree of the sub-hypergraph $H_c$ defined above. We see that this is the same as the maximum degree of the auxiliary multi-hypergraph $G$. Based on this reduction, we see that Theorem \ref{list-edge-coloring-general} follows from Theorem \ref{max-degree-thm}, and that Theorem \ref{list-edge-coloring-asymptotic} follows from Theorem \ref{delcourt-postle}. 

We refer to the above auxiliary multi-hypergraph $G = \bigsqcup_{c} H_c$ together with the set of color classes $E_i = \{ (e_i, c) : c \in L(e_i) \}$, $i \in [n]$, as the \textit{list edge-cover multi-hypergraph} of the original multi-hypergraph $H$ with list assignment $L$. We state the following necessary conditions for being a list edge-cover multi-hypergraph.

\begin{prop} \label{list-edge-cover-prop}
If a multi-hypergraph $G$ with color classes $E_1, \ldots, E_n$ is a list edge-cover multi-hypergraph of some other multi-hypergraph, then
\begin{itemize}
	\item[(a)] For every connected component $C$ and color class $E_i$, there is at most one edge in $E(C) \cap E_i$.
	\item[(b)] For every two distinct connected components $C, C'$ and distinct color classes $E_i, E_j$, if there exist edges $e \in E(C) \cap E_i$, $f \in E(C) \cap E_j$, $e' \in E(C') \cap E_i$, $f' \in E(C') \cap E_j$, then $|e \cap f| = |e' \cap f'|$.
\end{itemize}
\end{prop}

Proposition \ref{list-edge-cover-prop} is an easy consequence of the following two facts: the connected components of $G$ are all sub-hypergraphs of some multi-hypergraph $H$, and the edges of $G$ lying in the same color class $E_i$ correspond to the same edge in $H$. These conditions for being list edge-cover multi-hypergraph are not always sufficient, but they are useful to have in mind when we prove Theorem \ref{galvin-color-degree} below.

\subsection{Proof of Theorem \ref{galvin-color-degree}} \label{section-proof-of-galvin}

Now we prove Theorem \ref{galvin-color-degree}. Based on the above reduction of list edge-colorings to full rainbow matchings, to prove Theorem \ref{galvin-color-degree} it suffices to prove the following theorem.

\begin{thm} \label{galvin-rainbow-matching}
For every $\Delta \ge 2$, there exist a bipartite graph $G$ and edge-color classes $E_0, \ldots, E_n$ that is the list edge-cover graph of some bipartite graph, such that $G$ has maximum degree $\Delta$, $|E_i| = \Delta$ for every $i$ and there is no full rainbow matching.
\end{thm}

Our construction for Theorem \ref{galvin-rainbow-matching} will be shown as follows. First we will use Lemma \ref{join-lemma} to construct a bipartite graph $G_0$ and edge-color classes $\mathcal{P}_0$ that is the list edge-cover graph of some bipartite graph, such that there is no full rainbow matching, and all but one of the color classes in $\mathcal{P}_0$ have size $\Delta$. Then we will take multiple copies of $G_0, \mathcal{P}_0$ and apply Lemma \ref{join-lemma} to get a desired construction for Theorem \ref{galvin-rainbow-matching}. 

\begin{figure}
	\centering
	\includegraphics[scale=0.9]{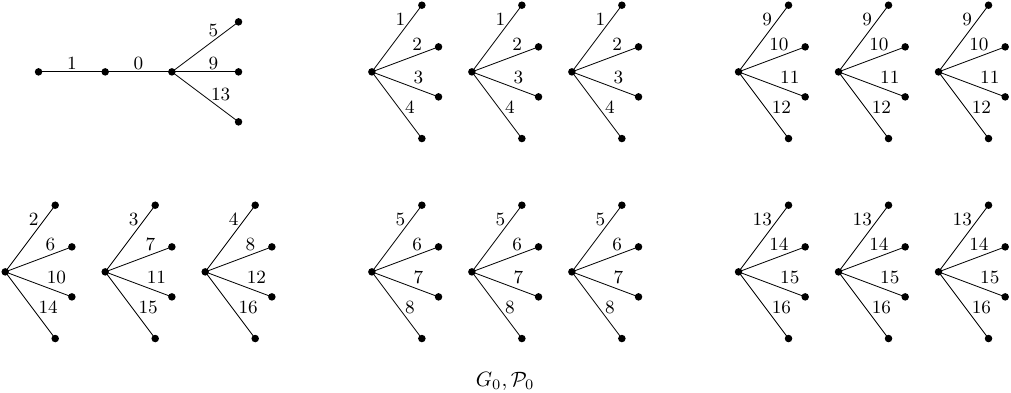}
	\caption{For maximum degree $\Delta = 4$, the list edge-cover graph $G_0$ and its edge-color classes $\mathcal{P}_0$ with no full rainbow matching. One color class has size 1, while the other color classes have size $\Delta$.}
	\label{list-edge-cover}
\end{figure}

\begin{figure}
	\centering
	\includegraphics[scale=1]{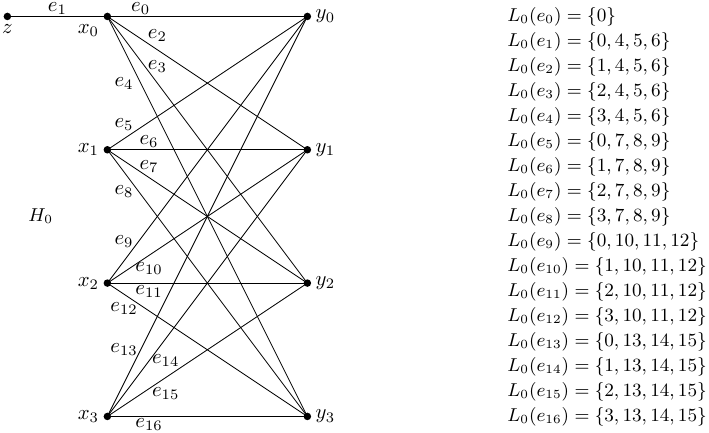}
	\caption{For maximum color degree $\Delta = 4$, the bipartite graph $H_0$ and list assignment $L_0$ for $E(H_0)$ with no proper $L_0$-coloring. One list has size $1$, while the other lists have size $\Delta$. The colors in the list $L_0(e_i)$ correspond to the connected components in Figure \ref{list-edge-cover} that contain an edge with color $i$.}
	\label{list-assignment}
\end{figure}

First we describe the bipartite graph $G_0$ and the color classes $\mathcal{P}_0 = \{F_0, F_1 \ldots, F_{\Delta^2}\}$, where $|F_0| = 1$ and $|F_i| = \Delta$ for $1 \le i \le \Delta^2$. The case $\Delta = 4$ is shown in Figure \ref{list-edge-cover}. The graph $G_0$ has two kinds of building blocks: one $\Delta$-broom, and many disjoint unions of stars $K_{1,\Delta}$. A \textit{$\Delta$-broom} $B$ is a tree with two central adjacent vertices $v$ and $w$, one leaf $u_0$ adjacent to $v$, and $\Delta - 1$ leaves $u_1, \ldots, u_{\Delta-1}$ adjacent to $v$. The graph $G_0$ has one copy of $B$, and we assign the edge $vw$ the color $0$, the edge $u_0v$ the color $1$, and the edge $u_iw$ the color $1+i \cdot \Delta$, for $i \in \{1, \ldots, \Delta-1\}$. Next, let forest $S$ be the disjoint union of $\Delta - 1$ copies of the star $K_{1,\Delta}$. In addition to $B$, the graph $G_0$ has $\Delta+1$ copies of $S$, say $S_0, \ldots, S_{\Delta}$. In the copy $S_0$, if its components are $K^1, \ldots, K^{\Delta-1}$, then the $\Delta$ edges in $K^i$ are given a unique color in $\{1 + i + j \cdot \Delta : j \in \{0, \ldots, \Delta - 1\}\}$. In the copy $S_i$, for $i \in \{1, \ldots, \Delta\}$, we properly edge-color each component with the colors $\{(i - 1) \cdot \Delta + j : j \in \{1, \ldots, \Delta\}\}$. In summary, the graph $G_0$ is the disjoint union of the subgraphs $B, S_0, S_1, \ldots, S_{\Delta}$ with the above edge-coloring. Letting $F_i$ be the set of edges that were given color $i$, for $i \in \{0, 1, \ldots, \Delta^2\}$, the resulting color classes $\mathcal{P}_0 = \{F_0, F_1 \ldots, F_{\Delta^2}\}$ have the sizes stated above.

We claim that $G_0$ has no full rainbow matching with respect to $\mathcal{P}_0$. This can be shown directly, but here we explain how $G_0, \mathcal{P}_0$ can be derived by iteratively applying Lemma \ref{join-lemma}, which helps to motivate their construction. We start with the $\Delta$-broom $B$ above and the edge-color classes $\mathcal{Q} = \{F_0, F_1''\}$, where $F_0 = \{vw\}$ and $F_1'' = E(B) \setminus \{vw\}$. Notice that $|F_0| = 1$, $|F_1''| = \Delta$, and that $B, \mathcal{Q}$ has no full rainbow matching. However, $B, \mathcal{Q}$ does not satisfy the necessary condition (a) in Proposition \ref{list-edge-cover-prop} for being a list edge-cover graph. To fix that, we proceed by iteratively adding copies of $S, \mathcal{P}'$ and applying Lemma \ref{join-lemma}, where $S$ is as above the disjoint union of $\Delta - 1$ copies of the star $K_{1, \Delta}$, and $\mathcal{P}'$ is the set of color classes of a proper edge-coloring of $S$ using $\Delta$ colors. (Observe that $S, \mathcal{P}'$ has no full rainbow matching, and that it is the list edge-cover graph of the star $K_{1, \Delta}$ with the uniform list assignment $L(e) = \{1, \ldots, \Delta - 1\}$, for $e \in E(K_{1, \Delta})$.) We add to $B, \mathcal{Q}$ one copy $S_0, \mathcal{P}_0'$ of the edge-colored graph $S, \mathcal{P}'$, and then we distribute each of the $\Delta$ edges of $F_1'' \in \mathcal{Q}$ into a unique class in $\mathcal{P}_0'$. The resulting edge-colored graph $B \sqcup S_0, \mathcal{Q}' = \{F_0, F_1',\ldots, F_\Delta'\}$ has no full rainbow matching, $|F_0| = 1$, $|F_i'| = \Delta$ for $i \in \{1, \ldots, \Delta\}$, and it satisfies the necessary condition (a) of Proposition \ref{list-edge-cover-prop} but not the necessary condition (b). Again, we fix that by adding copies of $S, \mathcal{P}'$. Iteratively for each color class $F_i'$, where $i \in \{1, \ldots, \Delta\}$, we add to the current edge-colored graph a copy $S_i, \mathcal{P}_i'$ of $S, \mathcal{P}'$, and then we distribute each of the $\Delta$ edges of $F_i'$ into a unique class in $\mathcal{P}_i'$. The result is the graph $G_0$ and color classes $\mathcal{P}_0 = \{F_0, F_1 \ldots, F_{\Delta^2}\}$ that we described above.

Next, we claim that $G_0, \mathcal{P}_0$ is the list edge-cover graph of some bipartite simple graph $H_0$ with list assignment $L_0$ for $E(H_0)$. An illustration of $(H_0, L_0)$ when $\Delta = 4$ is shown in Figure \ref{list-assignment}. The graph $H_0$ is the complete bipartite graph $K_{\Delta, \Delta}$, together with an external vertex and an edge connecting this vertex to some vertex in $K_{\Delta, \Delta}$. Labeling the parts of $K_{\Delta, \Delta}$ as $\{x_0, \ldots, x_{\Delta-1}\}$ and $\{y_0, \ldots, y_{\Delta-1}\}$, we label the edge $x_iy_j$ as $e_{1 + j + i \cdot \Delta}$ for $i, j \in \{0, \ldots, \Delta-1\}$, except for the edge $x_0y_0$ which we label as $e_0$. We connect the vertex $x_0$ to an external vertex $z$ and label the edge $x_0z$ as $e_1$. The list assignment $L_0$ for $E(H_0)$ is obtained by assigning a unique number to each connected component of $G_0$, and letting the list $L_0(e_i)$ consist of those numbers whose corresponding connected component has an edge of color $i$, for $i \in \{0, 1, \ldots, \Delta^2\}$. We then see that $G_0, \mathcal{P}_0$ is indeed the list edge-cover graph of $(H_0, L_0)$: The broom $B$ in $G_0$ comes from the edges $\{e_0\} \cup \{e_{1 + i \cdot \Delta} : i \in \{0, \ldots, \Delta-1 \}\}$, the collection of stars $S_0$ in $G_0$ arises from the collection of stars in $H_0$ centered at $y_1, \ldots, y_{\Delta-1}$, and each of the collection of stars $S_i$ in $G_0$ arises from a star in $H_0$ centered at $x_{i-1}$, for $i \in \{1, \ldots, \Delta\}$.

To complete our construction for Theorem \ref{galvin-rainbow-matching}, we only need to enlarge the color class $F_0$ of size $1$ into a class of size at least $\Delta$. We achieve this by iteratively applying Lemma \ref{join-lemma} to $G_0$. More directly, we take the disjoint union of $\Delta$ copies of $G_0, \mathcal{P}_0$ and then putting the $\Delta$ copies of $e_0$ into one common color class. Let $G$ be the resulting bipartite graph and let $\mathcal{P} = \{E_0, \ldots, E_{n}\}$ be the resulting color classes, where $n = \Delta^3$. Then we have that $|E_i| \ge \Delta$ for every $i$ and there is no full rainbow matching. Moreover, $G, \mathcal{P}$ is still the list edge-cover graph of some bipartite graph $H$ and list assignment $L$. Specifically, $(H, L)$ is obtained by taking $\Delta$ disjoint copies of $(H_0, L_0)$, then pasting together the disjoint copies of $H$ along the copies of $e_0$ (identifying all the copies of $x_0$, identifying all the copies of $y_0$, and then collapsing all the copies of $e_0$ into a single edge that we still call $e_0$), and letting $L(e_0)$ consist of the $\Delta$ distinct colors from the copies of list $L_0(e_0)$. Note that the resulting graph $H$ is still bipartite. Therefore, the edge-colored graph $G, \mathcal{P}$ gives the required construction for Theorem \ref{galvin-rainbow-matching}, and the associated bipartite graph $H$ and list assignment $L$ gives the required construction for Theorem \ref{galvin-color-degree}. This finishes the proof.

\bigskip

The above work on list edge-coloring is motivated by similar work on list vertex-coloring. In the vertex setting, Reed \cite{Re} conjectured that if a graph $H$ and list assignment $L$ for $V(H)$ has maximum (vertex) color degree $\Delta$, and $|L(v)| \ge \Delta+1$ for every vertex $v$, then there exists a proper $L$-coloring. This was suggested as a color degree generalization of the greedy upper bound on the list chromatic number, $\chi_{\ell}(H) \le \Delta(H)+1$. Haxell's theorem \cite{Ha2} on independent transversals implies that a proper $L$-coloring exists if $|L(v)| \ge 2\Delta$ for every vertex $v$, and Reed and Sudakov \cite{ReSu} proved that a proper $L$-coloring exists if $|L(v)| \ge (1 + o(1))\Delta$ for every vertex $v$ (see also \cite{LoSu}). However, Bohman and Holzman \cite{BoHo} found examples where $|L(v)| \ge \Delta+1$ for every $v$ and there is no proper $L$-coloring, which disproved Reed's conjecture. In \cite{HaWd1}, smaller counterexamples to Reed's conjecture were constructed in similar fashion to the above proof of Theorem \ref{galvin-color-degree}.

\section{Questions} \label{section-questions}
In this paper, we described a widely applicable method for constructing edge-colored multi-hypergraphs with large color classes and no full rainbow matchings. We focused on conditions depending on the maximum degree $\Delta$, and for example we were able to construct many edge-colored $r$-graphs with no full rainbow matchings achieving the extremal class sizes $r\Delta - 1$. But when it comes to rainbow matchings in properly edge-colored multigraphs, our investigations make it appear more relevant to look at the chromatic index instead of the maximum degree. We ask the following question.

\begin{ques} \label{rainbow-matching-question}
Let $G$ be a multigraph with chromatic index $\chi'$, and let $E_1, \ldots, E_n$ be the color classes of a proper edge-coloring of $G$. Is it true that if $|E_i| \ge (1 + o(1))\chi'$ for every $i$, then there exists a full rainbow matching? Could the lower bound $(1+o(1))\chi'$ be replaced by $\chi'+C$ for some constant $C$?
\end{ques}

Note that Theorem \ref{delcourt-postle} of Delcourt and Postle \cite{DePo} implies that the asymptotic lower bound is true if we assume that the multigraph has edge-multiplicity $o(\Delta)$ (in which case $\chi' = (1+o(1))\Delta$). Our construction for Theorem \ref{proper-edge-coloring-multigraph-1}, statement (3), implies that Question \ref{rainbow-matching-question} is not true with the constant $C = 2$, but we have not ruled out any $C \ge 3$. If we only consider bipartite graphs, then Theorem \ref{proper-edge-coloring-multigraph-1}, statement (2), implies that Question \ref{rainbow-matching-question} is not true with the constant $C = 1$, but we have not ruled out any value $C \ge 2$.

Positive answers to Question \ref{rainbow-matching-question} would lead to positive answers to the following question, which is about how the List Edge-Coloring Conjecture could extend to the ``color" setting. For a multigraph $H$, list assignment $L$ for $E(H)$, and a color $c$, recall that the subgraph $H_c$ has edge set $E_c = \{e \in E(H) : c \in L(e)\}$. The \textit{maximum color chromatic index} of $(H, L)$ is the maximum, over all colors $c$, of the chromatic index of $H_c$.

\begin{ques} \label{list-coloring-question}
Let $H$ be a multigraph, and let $L$ be a list assignment for $E(H)$. Is it true that if $(H, L)$ has maximum color chromatic index $\chi'$ and $|L(e)| \ge (1+o(1))\chi'$ for every $e \in E(H)$, then there exists a proper $L$-coloring of $H$? Could the lower bound $(1+o(1))\chi'$ be replaced by $\chi' + C$ for some constant $C$?
\end{ques}

Note that Theorem \ref{galvin-color-degree} implies that Question \ref{list-coloring-question} is not true with constant $C = 0$ even for bipartite graphs, but we have not ruled out any $C \ge 1$. Answers to Question \ref{rainbow-matching-question} and Question \ref{list-coloring-question} for bipartite multigraphs $H$ may be of particular interest. These questions could also be asked for general uniform multi-hypergraphs.

\section*{Acknowledgments}
This research was supported in part by the Austrian Science Fund (FWF) [10.55776/F1002]. This research was primarily done while the author was at the University of Waterloo in Waterloo, Canada. The author would like to thank the anonymous referees for their helpful comments and corrections.

\end{document}